\documentclass[a4paper,10pt]{amsart}
\usepackage{amssymb, amsmath, latexsym, amscd, verbatim, amsthm, color, tipa, enumerate, array}

\newtheorem{theorem}{Theorem}[section]
\newtheorem{lemma}[theorem]{Lemma}

\newtheorem{corollary}[theorem]{Corollary}
\newtheorem{proposition}[theorem]{Proposition}

\newtheorem*{thm}{Theorem}

\theoremstyle{definition}

\newtheorem{definition}[theorem]{Definition}

\newtheorem{example}[theorem]{Example}

\numberwithin{equation}{section}

\newcommand{\rank}{\textnormal{rank}}

\newcommand{\tr}{\textnormal{Tr}}

\newcommand{\ord}{\textnormal{ord}}

\newcommand{\set}[1]{\{\,#1\,\}}

\newcommand{\ifft}{if and only if }

\newcommand{\przystaje}[3]{#1 \equiv #2 \, \, (\mspace{1mu} \textnormal{mod} \mspace{1mu} \, #3)}

\newcommand{\cyclo}[1][]{\mathbb{Q}(\zeta_{#1})}

\newcommand{\spanp}{\textnormal{span}}

\newcommand{\bz}{\mathbb{Z}}
\newcommand{\bq}{\mathbb{Q}}
\newcommand{\br}{\mathbb{R}}

\newcommand{\fo}{\mathcal{O}}

\newcommand{\gtm}{\mathfrak{m}}

\newcommand{\fok}{\fo_{K}}

\newcommand{\trkq}{\tr_{K/\bq}}
\newcommand{\trcq}{\tr_{\cyclo /\bq}}
\newcommand{\trck}{\tr_{\cyclo / K}}

\title[Upper bounds for the Euclidean minima of abelian fields]{Upper bounds for the Euclidean minima of abelian fields of odd prime power conductor}
\author{Eva Bayer-Fluckiger, Piotr Maciak}
\address{Ecole Polytechnique F\'{e}d\'{e}rale de Lausanne\\EPFL--FSB--MATHGEOM--CSAG \\
1015 Lausanne\\
Switzerland}
\email{eva.bayer@epfl.ch}
\email{piotr.maciak@epfl.ch}
\date{\today}

\begin{document}

\begin{abstract}
The aim of this paper is to give upper bounds
for the Euclidean minima of abelian fields of
odd prime power conductor. In particular, these
bounds imply Minkowski's conjecture for totally
real number fields of conductor $p^r$, where $p$
is an odd prime number and $r \ge 2$.
\end{abstract}
\maketitle

\section{Introduction}
Let $K$ be an algebraic number field, and let $\fok$ be its ring of integers. Let ${\rm N}: K \to \bq$ be
the absolute value of the norm map. The number field $K$ is said to be {\it Euclidean} (with respect to the norm) if for every $a,b \in \fok$ with $b \not = 0$ there exist
$c, d \in \fok$ such that $a = bc + d$ and ${\rm N}(d) < {\rm N}(b)$.
It is easy to check that $K$ is Euclidean if and only if for every $x \in K$ there exists $c \in \fok$ such that ${\rm N}(x-c) < 1$. This suggests to look at
$$
M(K) = {\rm sup}_{x \in K} {\rm inf}_{c \in \fok} {\rm N}(x-c),
$$
called the {\it Euclidean minimum} of $K$.

The study of Euclidean number fields and Euclidean minima is a classical one,
see for instance~\cite{lemmermeyer}  for a survey. The present paper is
concerned with {\it upper bounds} for $M(K)$ in the case where $K$ is
an abelian field of odd prime power conductor. Let us recall some previous results.
Let $n$ be the degree of $K$ and $D_K$ the absolute value of its discriminant.
It is shown in~\cite{bayer} that for any number field $K$, we have $M(K) \le 2^{-n}D_K$.
The case of {\it totally real} fields is especially interesting, and has been the
subject matter of several papers. In particular, a conjecture attributed to Minkowski
states that if $K$ is totally real, then $$M(K) \le 2^{-n} \sqrt {D_K}.$$ This conjecture
is proved for $n \le 8$, cf.~\cite{grs7, grs8, mcmullen}; see also McMullen's paper~\cite{mcmullen} for a proof of the
case $n \le 6$ and a survey of the topic.

The point of view taken in the present paper is to study this conjecture for totally
real abelian fields, and more generally give upper bounds for Euclidean minima
of abelian fields. A starting point of this investigation is~\cite{bayer} where it is
proved that we have $M(K) \le 2^{-n} \sqrt {D_K}$ if $K$ is a cyclotomic field
of prime power conductor or the maximal totally real subfield of such a field.
The present paper contains some results concerning abelian fields
of odd prime power conductor. In particular, we show that if $K$ is such a field, then there exist constants $C=C(K) \leq \frac{1}{3}$ and
$\varepsilon = \varepsilon(K) \leq 2$ such that
\[
  M(K) \leq  C^n \, (\sqrt{D_K})^{\varepsilon}.
\]
If $[K :  \bq] > 2$, then one may choose $\varepsilon(K) < 2$.  Moreover, we show that $\varepsilon$ is asymptotically equal to $1$
and that under certain assumptions $C$ is asymptotically equal to $\frac{1}{2\sqrt{3}}$; see Theorem~(\ref{T:bound1}) for the precise statement.
In Theorem~(\ref{T:bound2}) we obtain the bound
\[
 M(K) \leq \omega^n \sqrt{D_K},
\]
where $\omega = \omega(K)$ is a constant which under certain assumptions is asymptotically equal to $\frac{1}{2\sqrt{3}}$.
In particular, using these bounds we show
\begin{thm}
 Suppose that K is a totally real field of conductor $p^r$, where $p$ is an odd prime and $r \ge 2$. Let $n$ be the degree of $K$ and let
$D_K$ be its discriminant. Then $$M(K) \le 2^{-n} \sqrt{D_K}.$$
\end{thm}
\noindent In other words, Minkowski's conjecture holds for this family of fields.

\bigskip
The strategy of the proofs is the following. If $K$ is an algebraic number field, we consider
 lattices defined on the ring of integers $\fok$ in the sense of \cite{ideal}, \S 1. This leads to a
 Hermite--like invariant of $\fok$, denoted by $\tau_{\rm min}(\fok)$, cf. \cite{ideal}, Definition 9, and
 \S 4 of the present paper. By
 \cite{bayer}, Corollary (5.2), we have
 $$M(K) \le  \left( {\tau_{\rm
min}(\fok) \over n } \right) ^{n/2} \sqrt{D_K},$$
where $n$ is the degree of $K$ and $D_K$ the absolute value of its discriminant.
In order to apply this result, we have to estimate $\tau_{\rm min}(\fok)$. The main technical
task of this paper is to do this in the case of abelian fields of odd prime power conductor.

\bigskip
The paper is structured as follows. After a brief section containing the notation used
throughout the paper, \S 3 describes the main results (theorems (\ref{T:bound1}) and (\ref{T:bound2})).
The rest of the paper is devoted
to their proofs, starting in \S 4 with a summary of some notions and results
concerning lattices and number fields, and their relation to Euclidean minima.
Suppose now that $K$ is an abelian field of prime power conductor. In \S 5,
we construct integral bases of $\fok$.
In \S 6, these bases are used to describe the lattice obtained by the canonical embedding
of $\fok$ (equivalently, the lattice given by the trace form). It turns out that this lattice is
isomorphic to the
orthogonal sum of lattices similar to the dual of a root lattice of type A and of a
lattice invariant by a symmetric group which already appears in [2]. Using this information, we obtain
an estimate of the Hermite--like thickness of the lattice $\fok$,  leading to an upper bound of $\tau_{\rm min}(\fok)$ that we apply in \S 7 to prove theorems~(\ref{T:bound1}) and (\ref{T:bound2}).
Finally, \S 8 contains some partial results and open questions concerning abelian fields
of odd prime conductor.

\section{Notation and a definition}
The following notation will be used throughout this paper. The set of all abelian extensions of $\bq$ of odd prime power conductor will be denoted by $\mathcal{A}$.
For $K \in \mathcal{A}$ we set:
\begin{align*}
 n &- \textnormal{ the degree of } K / \bq,\\
 D &- \textnormal{ the absolute value of the discriminant of } K,\\
 p &- \textnormal{ the unique prime dividing the conductor of } K, \\
 r &- \textnormal{ the } \textnormal{$p$-adic additive valuation of the conductor of } K,\\
 \zeta &- \textnormal{ a primitive root of unity of order } p^r,\\
 e &- \textnormal{ the degree  $[\cyclo : K]$}.
\end{align*}
If the dependence on the field $K$ needs to be emphasized, we shall add the index $K$ to the above symbols. For example, we shall write $n_K$ instead of $n$.
\bigskip
We also need the following definition

\begin{definition}
Let $\psi : \mathcal{D} \to \mathbb{R}$ be a function, where $\mathcal{D} \subset \mathcal{A}$. We shall say that
$\psi_o \in \mathbb{R}$ is the limit of $\psi$ as $n_K$
goes to infinity and write
\[
 \lim_{n_K \to \infty} \psi(K)  = \psi_0
\]
if for every $\epsilon > 0$ there exists $N >0$ such that for every field $K \in \mathcal{D}$
\[
 n_K > N  \implies  |\psi(K) - \psi_0| < \epsilon.
\]
We shall also write
\[
 \lim_{p_K \to \infty} \psi(K)  = \psi_0
\]
if for every $\epsilon > 0$ there exists $N >0$ such that for every field $K \in \mathcal{D}$
\[
p_K > N \implies |\psi(K) - \psi_0| < \epsilon.
\]
\end{definition}

\section{Euclidean minima -- statement of the main results}
In this section we present the main results of the paper; the proofs will be given in \S 7.
We keep the notation and definitions of the previous sections.

\begin{theorem}\label{T:bound1}
Let $K \in \mathcal{A}$. Then there exist constants $\varepsilon=\varepsilon(K) \leq 2$ and $C=C(K) \leq \frac{1}{3}$ such that
\[
  M(K) \leq  C^n \, (\sqrt{D_K})^{\varepsilon}.
\]
If $[K :  \bq] > 2$, then one may choose $\varepsilon(K) < 2$. Moreover,
\[
 \lim_{n_K \to \infty} \varepsilon(K) = 1.
\]
If $r_K \geq 2$, or $r_K=1$ and $[\cyclo:K]$ is constant, then we also have
\[
 \lim_{p_K \to  \infty} C(K) = \frac{1}{2\sqrt{3}}.
\]
\end{theorem}

\bigskip

\begin{theorem}\label{T:bound2}
Let $K \in \mathcal{A}$. Then there is a constant $\omega=\omega(K)$ such that
\[
 M(K) \leq \omega^n \sqrt{D_K}.
\]
If $r_K \geq 2$, or $r_K=1$ and $[\cyclo:K]$ is constant, then
\[
 \lim_{p_K \to \infty} \omega(K) = \frac{1}{2\sqrt{3}}.
\]
Moreover, if $r_K \geq 2$, then $\omega(K) \leq 3^{-2/3}$.
\end{theorem}

\bigskip

Note that this implies that Minkowski's conjecture holds for all totally real fields
 $K \in \mathcal{A}$ with composite conductor

 \begin{corollary}\label{T:Minkowski}

 Let $K \in \mathcal{A}$, and suppose that the conductor of $K$ is of the form $p^r$ with
 $r  > 1$. Then
\[
 M(K) \leq 2^{-n} \sqrt{D_K}.
\]

\end{corollary}

\bigskip
This follows from Theorem~(\ref{T:bound2}), since $3^{-2/3} < 1/2$, and for $K$ totally real this is precisely
Minkowski's conjecture.

\section{Lattices and number fields}

We start by recalling some standard notion concerning Euclidean lattices (see for instance
\cite {conway} and \cite {martinet}. A {\it lattice} is a pair $(L,q)$, where $L$ is a free $\bz$--module of finite rank, and  $q : L_{\br} \times L_{\br} \to \br$
is a positive definite symmetric bilinear form, where $L_{\br} = L \otimes_\bz \br$. 
If $(L,q)$ is a lattice and $a \in \br$, then we denote by $a(L,q)$ the lattice
$(L,aq)$. Two lattices $(L,q)$ and $(L',q')$ are said to be
{\it similar} if and only if there exists $a \in \br$ such that $(L',q')$ and $a(L,q)$ are isomorphic,
in other words if  there exists an isomorphism of $\bz$-modules $f : L \to L'$ such that $q'(f(x),f(y)) = a q(x,y)$.

\medskip

Let $(L,q)$ be a lattice, and set $q(x) = q(x,x)$. The {\it maximum }of $(L,q)$ is defined by
\[
 {\rm max}(L,q) = \sup_{x \in L_{\br}} \inf_{c \in L} q(x-c).
\]
Note that ${\rm max}(L,q)$ is the square of the covering radius of the associated sphere covering.
The {\it determinant} of $(L,q)$ is denoted by ${\rm det}(L,q)$. It is by definition the determinant of the matrix of $q$ in a $\bz$--basis of $L$.
The {\it Hermite--like thickness} of $(L,q)$ is
$$
\tau(L,q) = {{\rm max}(L,q) \over {\rm det}(L,q)^{1/m}},
$$
where $m$ is the rank of $L$. Note that $\tau(L,q)$ only depends on the similarity class of the lattice $(L,q)$.

\bigskip
{\it A family of lattices}
\medskip

Let $m\in {\mathbb {N}}$, and $b\in {\br}$ with $b>m$.
Let $L=L_{b,m}$ be a lattice in ${\br} ^m$ with Gram matrix
$$b I_m - J_m = \left( \begin{array}{cccc}
b-1 & -1 & \ldots  & -1 \\
-1 & \ddots & \ddots & \vdots  \\
\vdots & \ddots & \ddots & -1  \\
-1 & \ldots & -1 & b-1
\end{array} \right),$$
where $I_m$ is the $m\times m$-identity matrix and $J_m \in \{1 \} ^{m\times m}$ is the all-ones matrix.
Then $L$ is a lattice of determinant  $(b-m)b^{m-1}$. Moreover the automorphism group of $L$ contains $\langle -I_m \rangle \times S_m$, where
the symmetric group $S_m$ acts by permuting the coordinates. These lattices were defined in \cite{bayer-nebe}, (4.1).
Note that the lattice $L_{m+1,m} $ is similar to the dual lattice $A_m^{\#}$ of the
root lattice $A_m$ (see for instance \cite{conway}, Chapter 4, \S 6, or \cite{martinet} for the definition of the root lattice $A_m$).

\bigskip
{\it Lattices defined over number fields}

\medskip In the sequel, we will be concerned with lattices defined on rings of integers
of abelian number fields.
Let $K$ be an number field of degree $n$, and suppose that $K$ is either totally
real or totally complex. Let us denote by
$^{\overline {\ }} : K  \to K$ the identity in the first case and the
complex conjugation in the second one, and let $P$ be the set of totally
positive elements of the fixed field of this involution.  Let us denote by ${\rm Tr} : K \to \bq$ the trace
map. For any $\alpha \in P$, set $q_{\alpha}(x,y) = {\rm Tr}(\alpha x \overline y)$ for all $x, y \in K$.
Then $(\fok,q_{\alpha})$ is a lattice.  Set
$$
\tau_{\rm min}(\fok) = {\rm inf} \{ \tau(\fok,q_{\alpha}) \ | \ \alpha \in P \}.
$$

If $D_K$ is the absolute value of the discriminant of $K$, then, by \cite{bayer}, Corollary (5.2), we have
\begin{equation}\label{E:gen-est}
 M(K) \leq \left(\frac{\tau_{\min}(\fok)}{n}\right)^{\frac{n}{2}} \sqrt{D_K},
\end{equation}
This upper bound will be used in \S 7 to prove theorems~(\ref{T:bound1}) and (\ref{T:bound2}).
\section{Gaussian periods and integral bases}

Let $K \in \mathcal{A}$. In order to exploit the upper bound of \S 4, we need some information
concerning the lattices defined on the ring of integers $\fok$, and these will be described
using integral bases of $\fok$. The aim of this section is to find such bases. This will be done
in the spirit of the work of Leopold~\cite{leopold}, see also Lettl~\cite{lettl}.

\medskip
Recall that the $p^r$ is the conductor of $K$  and that $e=[\cyclo:K]$.  Then
$e$ divides $p-1$. This implies that 
the extension $\cyclo / K$ is tamely ramified, and hence 
the trace map ${\rm Tr} : \bz[\zeta] \to \fo_K$ is surjective.

\medskip
Set $R=\bz / p^{r} \bz$ and 
let us denote by $H$ the unique subgroup of order $e$ of $R^{\ast}$. Then $H$ acts
on $R$ by left multiplication. The orbit $H 0$ will be denoted by $\mathbf{0}$ and called the zero orbit.

\begin{definition}
For $\mathbf{x} \in R /H \setminus \set{\mathbf{0}}$,  we define a \textit{Gaussian period}
\[
 f_{\mathbf{x}} = \sum_{x \in \mathbf{x}} \zeta^x.
\]
In addition, we set $f_{\mathbf{0}} = e$.
\end{definition}

If $\mathbf{x}=H x$, then $f_{\mathbf{x}} = \tr_{\cyclo / K}(\zeta^x)$. As  the trace map ${\rm Tr} : \bz[\zeta] \to \fo_K$ is surjective., Gaussian periods are generators of $\fok$ over $\bz$. The next
proposition will be used to show that $\fok$ has actually an integral basis consisting of Gaussian periods.  

\medskip
Set  $S=\bz / p^{r-1} \bz$, and let
$\pi : R \to S$ be the canonical projection. The group $H$ acts on $S$ by $h \cdot s = \pi(h)  s$.
Clearly, $\pi$ is a morphism of $H$-sets and hence it induces the unique map between the orbit sets $\rho : R / H \to S / H$ such that $\rho(H x) = H \pi(x)$ for all $x \in R$.
In other words, if $\mu_R : R \to R / H$ and $\mu_S : S \to S / H$ are the canonical projections, then $\rho \mu_R = \mu_S \pi$. In particular, $\rho$ is surjective.
For any subset $A$ of $R$, set $\zeta^A = \{\zeta^a \ | \  a \in A \}$, and let us denote by $A^c$ the
complement of $A$ in $R$. For a finite set $X$, we denote by $|X|$ the number of elements of $X$.

\medskip
 We thank H.W. Lenstra, Jr, for sending us the part  (1) $\Leftrightarrow$ (2)  of
the following proposition.

\begin{proposition}\label{P:set-A}
Let $A \subset R \setminus \set{0}$ be an $H$-invariant set. The following conditions are equivalent:
\begin{enumerate}
 \item $\zeta^A$ is a basis of $\bz[\zeta]$.
 \item The restriction $\pi :A^c \to S$ is a bijection.
 \item The restriction $\rho : A^c / H \to S / H$ is a bijection.
 \item\label{P:set-A-card} For every $\mathbf{y} \in S / H$ we have
 \[
   |\rho^{-1}(\mathbf{y}) \cap A / H| = |\rho^{-1}(\mathbf{y})| -1.
 \]
 \item $\trck(\zeta^{A}) = \set{f_{\mathbf{x}} \mid \mathbf{x} \in A/ H}$ is a basis of $\fok$.
\end{enumerate}

\end{proposition}

\noindent
{\it Proof}.
 (1) $\Leftrightarrow$ (2) Note that the sum of powers of $\zeta$ over any coset of $\ker \pi$ in $R$ equals zero. Thus, if $\zeta^A$ is a $\bz$-basis, then $A$ must miss at least one element of each coset.
It cannot miss more than one since then the cardinality of $A$ would be too small. Conversely, if $A$ misses exactly one element from
each coset, then the sum relation mentioned shows that the $\bz$-span of $\zeta^A$ contains all roots of unity of order $p^r$. Hence $\zeta^A$ forms an integral basis. 

\medskip

\noindent (2) $\Rightarrow$ (3) This follows immediately from the fact that $\pi$ is a morphism of $H$-sets and $A^c$ is $H$-invariant.

\medskip

\noindent (3) $\Rightarrow$ (2) First, we shall show that $\pi_{|A^c}$ is onto. Suppose that this is not true. Then $\pi$ must miss at least one full orbit since it is an $H$-map
but in such a case the restriction of $\rho: A^c / H \to S / H$ would not be surjective. Thus $\pi_{|A^c}$ is onto. We claim that $|A^c| = |S|$, which implies that $\pi_{|A^c}$
is a bijection. Indeed, the set $A^c /H$ maps bijectively onto $S/H$, which implies that they have the same cardinality. Since both $A^c /H$ and $S /H$ contain the respective zero orbits, it is easy
to check that
\[
 \left|A^c / H \right| = 1 + \frac{|A^c|-1}{e} \quad \textnormal{ and } \quad \left|S / H \right| = 1 + \frac{|S|-1}{e},
\]
which readily implies $|A^c| = |S|$.

\medskip

\noindent (3) $\Leftrightarrow$ (4) Since $A$ is $H$-invariant, it follows that the sets $A/H$, $A^c /H$ form a partition of the orbit space $R /H$. Consequently, for every $\mathbf{y} \in S / H$ we have
\[
 |\rho^{-1}(\mathbf{y}) \cap A / H| = |\rho^{-1}(\mathbf{y})| - |\rho^{-1}(\mathbf{y}) \cap A^c / H|.
\]
The restriction $\rho : A^c / H \to S / H$ is a bijection \ifft $|\rho^{-1}(\mathbf{y}) \cap A^c / H| = 1$ for every $\mathbf{y} \in S / H$.

\medskip

\noindent (1) $\Rightarrow$ (5) Assume now that $\zeta^A$ is a basis of $\bz[\zeta]$. We shall show that  $\trck(\zeta^{A})$ is an integral basis of $K$. Since $0 \notin A$, it follows that $A$ is a union of $n$ orbits, each of cardinality $e$.
Consequently,
\[
 |\trck(\zeta^{A})| \leq n = \rank \, \fok.
\]
Since $\cyclo / K$ is tamely ramified, $\trck(\zeta^A)$ generates $\fok$, which in turn implies that we have in fact $|\trck(\zeta^{A})|  = n$
and that $\trck(\zeta^{A})$ is an integral basis of $K$. 

\medskip

\noindent (5) $\Rightarrow$ (1) Since $\trck(\zeta^A)$ is a basis, we have $|\trck(\zeta^{A})|  = n$. It follows that $|A| \geq n e = p^{r-1} (p-1)$. Suppose by contradiction that $\zeta^A$ is not a basis.
Then $\pi : A^c \to S$ is not a bijection. Since $$|A^c| = |R| - |A| \leq p^r - p^{r-1} (p-1) = p^{r-1} = |S|,$$ it follows that $\pi : A^c \to S$ is not a surjection. Consequently,
there exists $x_0 \notin \ker \pi$ such that the coset $C=x_0 + \ker \pi$ is contained in $A$. Note that the intersection of $C$ with any $H$-orbit is either empty or a singleton. Indeed, if
$x_0 + z_1 = h (x_0 + z_2)$ for some $z_1, z_2 \in \ker \pi$ and $h \in H \setminus \set{1}$, then $(1-h) x_0 \in \ker \pi$. Note that $1-h$ is invertible. Hence $x_0 \in \ker \pi$,
which is a contradiction. Therefore the set $HC = \set{hc \mid h \in H, c \in C}$ is contained in $A$ and it has $|H| \cdot |C| $ elements. Furthermore,
$HC /H$ is contained in $A / H$ and
\[
 \sum_{\mathbf{x} \in HC /H} f_{\mathbf{x}} = \sum_{hc \in HC} \zeta^{hc} = \sum_{h \in H} \sum_{c \in C} \zeta^{hc} =  0,
\]
which contradicts the linear independence of elements of $\trck(\zeta^A)$. Thus $\zeta^A$ is a basis of $\bz[\zeta]$.

\bigskip
This proposition implies that $\fok$ has an integral basis consisting of Gaussian periods.
Indeed, we have

\begin{corollary}\label{C:basis}
 There exists an $H$--invariant set $A \subset R \setminus \set{0}$
 such that
$$
\trck(\zeta^{A}) = \set{f_{\mathbf{x}} \mid \mathbf{x} \in A/ H}
$$
 is a basis of $\fok$.
  \end{corollary}

\noindent
{\it Proof.} For all $\mathbf{y} \in S/ H$ with $\mathbf{y} \not = \mathbf{0}$, let
us choose $\mathbf{x}_{\mathbf{y} } \in R/ H$ such that $\rho(\mathbf{x}_{\mathbf{y} } ) =
\mathbf{y}$.
Set $\mathbf{x}_{\mathbf{0}} = \mathbf{0}$, and let $B = \cup_{\mathbf{y} \in S/ H}\mathbf{x}_{\mathbf{y} }$. Then $B$ is an $H$--invariant subset of $R$ containing $0$, and
the restriction $\rho : B / H \to S / H$ is a bijection. Set $A = B^c$; then $A$
is an $H$--invariant subset of $R \setminus \set{0}$, and the restriction
$\rho : A^c / H \to S / H$ is a bijection. By Proposition~(\ref{P:set-A}), this implies that
$\trck(\zeta^{A}) = \set{f_{\mathbf{x}} \mid \mathbf{x} \in A/ H}$ is a basis of $\fok$.

\bigskip
\bigskip

\section{Geometry of the ring of integers}

We keep the notation of the previous section; in particular, $K \in \mathcal{A}$ and
$p = p_K$. Recall that $\fok$ is the ring of integers of $K$, and let us consider the
lattice $(\fok,q)$, where $q$ is defined by $q(x,y) = \trkq(x \overline y)$. As we have
seen in \S 4, the Hermite--like thickness of this lattice can be used to give an upper bound of
the Euclidean minimum of $K$. The purpose
of this section is to describe the lattice $(\fok,q)$ using the results of \S 5, so that we can
compute its Hermite--like thickness.

\bigskip
We will see that $(\fok,q)$ decomposes in a natural way into the
orthogonal sum of a lattice $\Gamma_K$, which is similar to the orthogonal sum of copies of  the dual lattice of
the root lattice $A_{p-1}$, and of a lattice $\Lambda_K$,  which is similar to a certain lattice of type $L_{b,m}$ defined in \S 4.
The Hermite--like thickness of these lattices can be estimated, cf.~\cite{bayer-nebe}, Theorem (4.1).
This allows us to give good upper bounds for the Euclidean minima of fields $K \in \mathcal{A}$,
following the strategy outlined in the introduction and in \S 4.


\bigskip Let $\Gamma_K$  be the orthogonal sum of
$\frac{p^{r-1}-1}{e}$ copies of the lattice $p^{r-1}A_{p-1}^{\#}$.
Set  $d=\frac{p-1}{e}$, and let  $\Lambda_K = e p^{r-1} L_{\frac{p}{e},d}$ (note that
the scaling is taken in the sense of \S 4, that is it refers to multiplying the quadratic form by
the scaling factor). 

\begin{theorem}\label{T:ortho-sum}
The lattice $(\fok,q)$ is isometric to the orthogonal sum of $\Gamma_K$ and of $\Lambda_K$.
\end{theorem}

Before proving this theorem, we need a few lemmas.
Recall that $R$ denotes
the ring $\bz / p^r \bz$, and let $\gtm = p \bz /p^r \bz$ be the maximal ideal of $R$. Note that if $r=1$, then $\gtm$ is  the zero ideal.
For $\mathbf{x} \in R / H$ we set $\ord_p(\mathbf{x}) = {\rm max} \{ k \in {\bf N} \ |  \ \mathbf{x} \subset \gtm^k \}$. Let us denote by $\mu$ the M\"obius function.

\begin{lemma}\label{L:trace-of-gp}
Let $\mathbf{x} \in R / H$. Then,
 \[
  \trkq(f_{\mathbf{x}}) = \frac{\phi(p^r)}{\phi(p^{r-s})} \cdot \mu(p^{r-s}),
 \]
where $s = \ord_p(\mathbf{x})$
\end{lemma}

\noindent
{\it Proof.}
Let $x_0 \in \mathbf{x}$. We have
\[
  \trkq(f_{\mathbf{x}}) = \trkq(\trck(\zeta^{x_0})) = \trcq(\zeta^{x_0})
\]
 Assume first that $s=0$. Then, $\mathbf{x} \subset R^{\ast}$ and hence $x_0 \in R^{\ast}$. Consequently,
\[
 \trkq(f_{\mathbf{x}}) = \sum_{x \in R^{\ast}} \zeta^{x_0 x} = \sum_{x \in R^{\ast}} \zeta^{x} = \mu(p^r) = \frac{\phi(p^r)}{\phi(p^{r-s})} \cdot \mu(p^{r-s}).
\]
Now, assume that $1 \leq s < r$. Then, $x_0 = p^s x_1$ with $x_1 \in R^{\ast}$. Set $\xi = \zeta^{p^s}$ and $T=\bz / p^{r-s} \bz$. Then $\xi$ is a primitive root of unity of order $p^{r-s}$.
If $\tau : R^{\ast} \to T^{\ast}$ is the natural map with  kernel $G$, and the set $Y \subset R^{\ast}$ is mapped by $\tau$ bijectively onto $(\bz / p^{r-s} \bz)^{\ast}$, then
\begin{align*}
 \trkq(f_{\mathbf{x}}) &= \sum_{x \in R^{\ast}} \xi^{x_1 x} = \sum_{x \in R^{\ast}} \xi^x = \sum_{g \in G} \sum_{y \in Y} \xi^{gy} =  \sum_{g \in G} \sum_{y \in Y} \xi^{y}\\
		       &= |G| \cdot \sum_{t \in T^{\ast}} \xi^{t} = \frac{\phi(p^r)}{\phi(p^{r-s})} \cdot \mu(p^{r-s}).
\end{align*}
Finally, if $s=r$, then $\mathbf{x}=\mathbf{0}$ and $x_0=0$ and hence
\[
   \trkq(f_{\mathbf{x}}) = \trcq(1) = \phi(p^r) = \frac{\phi(p^r)}{\phi(p^{r-s})} \cdot \mu(p^{r-s}).
\]

\bigskip

\begin{proposition}\label{P:Gram-prep}
 Let $\mathbf{x_1}, \mathbf{x_2} \in R / H \setminus \set{\mathbf{0}}$. Then,
\[
 \trkq(f_{\mathbf{x_1}} \overline{f_{ \mathbf{x_2}}}) =
\begin{cases}
p^r-p^{r-1} \quad &\text{ if }\, \mathbf{x_1} = \mathbf{x_2} \text{ and } \rho(\mathbf{x_1})=\rho(\mathbf{x_2}) \neq \mathbf{0},\\
p^r - e p^{r-1} \quad &\text{ if }\, \mathbf{x_1} = \mathbf{x_2} \text{ and } \rho(\mathbf{x_1})=\rho(\mathbf{x_2}) = \mathbf{0},\\
-p^{r-1} \quad &\text{ if }\, \mathbf{x_1} \neq \mathbf{x_2} \text{ and } \rho(\mathbf{x_1})=\rho(\mathbf{x_2}) \neq \mathbf{0},\\
-e p^{r-1} \quad &\text{ if }\, \mathbf{x_1} \neq \mathbf{x_2} \text{ and } \rho(\mathbf{x_1})=\rho(\mathbf{x_2}) = \mathbf{0},\\
0 \quad &\text{ if }\, \mathbf{x_1} \neq \mathbf{x_2} \text{ and } \rho(\mathbf{x_1}) \neq \rho(\mathbf{x_2}).
\end{cases}
\]
\end{proposition}

\noindent
{\it Proof.}
 Let $x_1 \in \mathbf{x_1}$ and $x_2 \in \mathbf{x_2}$.  For $h \in H$ we set $\mathbf{x}(h) = H (x_1-x_2h)$ and $s(h) = \ord_p \mathbf{x}(h)$.
Then,
 \begin{align*}
  f_{\mathbf{x_1}} \overline{f_{\mathbf{x_2}}} &= (\sum_{h_1 \in \hbar} \zeta^{x_1h_1}) (\sum_{h_2 \in \hbar} \zeta^{-x_2h_2}) = \sum_{h_1 \in \hbar} \sum_{h_2 \in \hbar} \zeta^{x_1h_1-x_2h_2}
	  = \sum_{h_1 \in \hbar} \sum_{h \in \hbar} \zeta^{x_1h_1-x_2hh_1}\\ &= \sum_{h_1 \in \hbar} \sum_{h \in \hbar} \zeta^{(x_1-x_2h)h_1}
	  = \sum_{h \in \hbar} \sum_{h_1 \in \hbar} \zeta^{(x_1-x_2h)h_1} = \sum_{h \in \hbar} f_{\mathbf{x}(h)}.
 \end{align*}
By Lemma~(\ref{L:trace-of-gp}), we have
\begin{equation}\label{E:tr-form}
\trkq(f_{\mathbf{x_1}} \overline{f_{ \mathbf{x_2}}}) = \sum_{h \in H} \trkq(f_{\mathbf{x}(h)}) =  \sum_{h \in H} \frac{\phi(p^r)}{\phi(p^{r-s(h)})} \cdot \mu(p^{r-s(h)}).
\end{equation}
If $\mathbf{x_1} = \mathbf{x_2}$, we can take $x_1 = x_2$ and then $\mathbf{x}(h) = H x_1(1-h)$. Clearly, $s(1) = r$. If $h\neq 1$, then  we have
$s(h) = \ord_p\mathbf{x_1}$. Thus, if $\ord_p\mathbf{x_1} < r-1$, then $\rho(\mathbf{x_1}) \neq \mathbf{0}$ and the only non-zero term of the sum~(\ref{E:tr-form})
is the one corresponding to $h=1$. Thus, we have
\[
 \trkq(f_{\mathbf{x_1}} \overline{f_{ \mathbf{x_2}}}) = \frac{\phi(p^r)}{\phi(p^{r-s(1)})} \cdot \mu(p^{r-s(1)}) = \phi(p^r) = p^r-p^{r-1}.
\]
If $\ord_p(\mathbf{x_1}) = r-1$, then $\rho(\mathbf{x_1}) = \mathbf{0}$ and the sum~(\ref{E:tr-form}) becomes
\begin{align*}
  \trkq(f_{\mathbf{x_1}} \overline{f_{ \mathbf{x_2}}}) &= \sum_{h \in H} \frac{\phi(p^r)}{\phi(p^{r-s(h)})} \cdot \mu(p^{r-s(h)})\\
       & = \phi(p^r) + (e-1) \cdot \frac{\phi(p^r)}{\phi(p^{r-1})} \cdot \mu(p)\\
      &= p^r - e p^{r-1}.
\end{align*}
Suppose now that $\mathbf{x_1} \neq \mathbf{x_2}$. Observe that $\rho(\mathbf{x_1}) = \rho(\mathbf{x_2})$ \ifft there is an $h \in H$ such that $s(h) =r-1$. Moreover,
 in such a case an element $h$ with this property is unique unless $\rho(\mathbf{x_1}) = \mathbf{0}$,
in which case we have $s(h) = r-1$ for
all $h \in H$. Thus, assuming that $\rho(\mathbf{x_1}) = \rho(\mathbf{x_2})$ and $\rho(\mathbf{x_1}) \neq \mathbf{0}$, we have
\begin{align*}
  \trkq(f_{\mathbf{x_1}} \overline{f_{ \mathbf{x_2}}}) &= \sum_{h \in H} \frac{\phi(p^r)}{\phi(p^{r-s(h)})} \cdot \mu(p^{r-s(h)})\\
						       &=  \frac{\phi(p^r)}{\phi(p)} \cdot \mu(p) = -p^{r-1}.
\end{align*}
If $\rho(\mathbf{x_1}) = \mathbf{0}$, then
\begin{align*}
 \trkq(f_{\mathbf{x_1}} \overline{f_{ \mathbf{x_2}}}) &= \sum_{h \in H} \frac{\phi(p^r)}{\phi(p^{r-s(h)})} \cdot \mu(p^{r-s(h)})\\
						      &= \sum_{h \in H} \frac{\phi(p^r)}{\phi(p)} \cdot \mu(p) = -e p^{r-1}.
\end{align*}
Finally, if $\rho(\mathbf{x_1}) \neq \rho(\mathbf{x_2})$, then $s(h) \leq r-2$ for all $h \in H$, which gives
\[
  \trkq(f_{\mathbf{x_1}} \overline{f_{ \mathbf{x_2}}}) = 0.
\]
\bigskip

\begin{lemma}\label{L:rho-inv}
 Let $\mathbf{y} \in S / H$. Then,
\[
 |\rho^{-1}(\mathbf{y})| =
\begin{cases}
 1+\frac{p-1}{e} \quad  &\text{ if } \, \mathbf{y}=\mathbf{0},\\
 p \quad  &\text{ if } \, \mathbf{y} \neq \mathbf{0}.
\end{cases}
\]
\end{lemma}

\noindent
{\it Proof.}
Let $X = \rho^{-1}(\mathbf{y})$. Note that $\mathbf{0} \in X$ if and only if $\mathbf{y}=\mathbf{0}$,
hence
\begin{equation*}
 |\mu_R^{-1}(X)| =
\begin{cases}
 1 + e \cdot(|X|-1) \quad &\text{ if } \, \mathbf{y}=\mathbf{0},\\
 |X| \cdot e \quad &\text{ if } \, \mathbf{y} \neq \mathbf{0}.
\end{cases}
\end{equation*}
On the other hand,
$
\mu_R^{-1}(X) = \mu_R^{-1}(\rho^{-1}(\mathbf{y})) = (\rho \mu_R)^{-1}(\mathbf{y}) = (\mu_S \pi)^{-1}(\mathbf{y}).
$

If $\mathbf{y}=\mathbf{0}$, then
\[
 \mu_R^{-1}(X) = (\mu_S \pi)^{-1}(\mathbf{0}) = \set{x \in R \mid H\pi(x) = \mathbf{0}} = \gtm^{r-1}
\]
and hence $|(\mu_S \pi)^{-1}(\mathbf{0})| = p$, which implies that $|X| = 1+\frac{p-1}{e}$.

If $\mathbf{y} \neq \mathbf{0}$, then there is an element $x_0 \in R \setminus \gtm^{r-1}$  such that $\mathbf{y} = H \pi(x_0)$. Consequently,
\begin{align*}
 \mu_R^{-1}(X) &=  (\mu_S \pi)^{-1}(\mathbf{y}) = \set{x \in R \mid H\pi(x) = H \pi(x_0)}\\
		    &= \set{x \in R \mid x = k + hx_0 \text{ for some } k \in \gtm^{r-1} \text{ and } h \in H  }.
\end{align*}
Let $k_1, k_2 \in \gtm^{r-1}$ and $h_1, h_2 \in H$. If $k_1 + h_1 x_0 = k_2 + h_2 x_0$, then $h_1 (1 - h_1^{-1}h_2 ) x_0 \in  \gtm^{r-1}$. Since $x_0 \notin \gtm^{r-1}$, it follows that
$1 - h_1^{-1}h_2$ is not invertible. Hence  $h_1=h_2$, which in turn implies that $k_1=k_2$. 
Therefore
$|(\mu_S \pi)^{-1}(\mathbf{y})| = |\gtm^{r-1}| \cdot |H| = p e$, which gives $|X| = p$.
\bigskip

\noindent
{\it Proof of Theorem~(\ref{T:ortho-sum}).}
Let $A \subset R \setminus \set{0}$ be an $H$-invariant set such that $\trck(\zeta^{A})$ is a basis of $\fok$. For $\mathbf{y} \in S / H$ we set
\[
 B_{\mathbf{y}} = \set{f_{\mathbf{x}} \in \trck(\zeta^{A}) \mid \rho({\mathbf{x}}) = \mathbf{y}} \quad \textnormal{ and } \quad L_{\mathbf{y}} = \spanp_{\bz} B_{\mathbf{y}}.
\]
By Proposition~(\ref{P:Gram-prep}), we have
\begin{equation}\label{E:OK-decomp}
 \fok = \perp_{\mathbf{y} \in S /H} L_{\mathbf{y}},
\end{equation}
in other words the lattice $(\fok, q)$ is the orthogonal sum of the lattices obtained by the restriction of $q$ to $L_{\mathbf{y}}$ for $\mathbf{y} \in S/H$.
Combining Lemma~(\ref{L:rho-inv}) and the condition (\ref{P:set-A-card}) of Proposition~(\ref{P:set-A}), we obtain that
$$
|B_{\mathbf{y}}|=
\begin{cases}
 \frac{p-1}{e} \quad  &\text{ if } \, \mathbf{y}=\mathbf{0},\\
 p-1 \quad  &\text{ if } \, \mathbf{y} \neq \mathbf{0}.
\end{cases}
$$
Furthermore, using Proposition~(\ref{P:Gram-prep}) again, we conclude that the Gram matrix of the lattice $L_{\mathbf{y}}$ with respect to $B_{\mathbf{y}}$ is $p^{r-1} (p I_{p-1} - J_{p-1})$ unless
$\mathbf{y} = \mathbf{0}$ in which case it equals $e p^{r-1} (\frac{p}{e} I_d - J_d)$, where $d=\frac{p-1}{e}$. Consequently, we have $L_{\mathbf{0}} = \Lambda_K$. Moreover, $S/H$ has $\frac{p^{r-1}-1}{e}$ nonzero orbits. As a result,
\[
 \perp_{\mathbf{y} \neq \mathbf{0}} L_{\mathbf{y}} \simeq \Gamma_K.
\]
Thus the equality~(\ref{E:OK-decomp}) implies that $(\fok,q)$ is isometric to the orthogonal sum of $\Gamma_K$ and $\Lambda_K$.
\bigskip

We now apply Theorem~(\ref{T:ortho-sum}) to give an upper bound of the Hermite--like thickness of the lattice $(\fok,q)$. The following is well--known

\begin{lemma}\label{L:det-fok}
 We have
\[
 \det (\fok,q) = p^{\upsilon},
\]
where
\[
 \upsilon = rn -\frac{(p^{r-1}-1)}{e} -1.
\]
\end{lemma}

\noindent {\it Proof.} Note that $\det (\fok,q)$ is the absolute value of the discriminant of $K$. 
The result follows from Theorem (4.1) in~\cite{disc}. Alternatively, one can compute $\det (\fok,q)$ directly using Theorem~(\ref{T:ortho-sum}).

\begin{lemma}\label{L:max-fok}
 We have
\[
 \max(\fok,q) \leq n \cdot \frac{p^{r+1}+p^r+1-e^2}{12 p}.
\]
\end{lemma}

\noindent {\it Proof.}
 By Theorem (4.1) in~\cite{bayer-nebe}, we have
\[
 \max (L_{\frac{p}{e},d}) \leq \frac{d(p^2+p+1-e^2)}{12ep}.
\]
Furthermore, $\max (L_{p, p-1}) = \frac{p^2-1}{12}$. Consequently,
\begin{align*}
 \max (\fok,q) = \sum_{\mathbf{y} \in S / H} \max(L_{\mathbf{y}}) &\leq \frac{p^{r-1}-1}{e} \cdot p^{r-1} \cdot \frac{p^2-1}{12} + p^{r-1} \cdot \frac{d(p^2+p+1-e^2)}{12p}\\
	     &= d p^{r-1} \left[ \left(\frac{p^{r-1}-1}{p-1}\right) \cdot \left(\frac{p^2-1}{12}\right) + \frac{p^2+p+1-e^2}{12 p} \right]\\
	     &= n \cdot \frac{(p^{r-1}-1)(p+1)p+(p^2+p+1-e^2)}{12p}\\
	     &= n \cdot \frac{p^{r+1}+p^r+1-e^2}{12 p}.
\end{align*}

As a direct consequence of the above lemmas, we obtain the following 
upper bound
of $\tau_{\rm min}(\fok)$

\begin{corollary}\label{c:tau-min}
We have
\[
 \tau_{\rm min}(\fok) \leq  \tau(\fok,q)   \leq n \cdot p^{r-\frac{\upsilon}{n}} \cdot \frac{p^{r+1}+p^r+1-e^2}{12 p^{r+1}}.
\]

\end{corollary}

\bigskip
This bound will be used in the next section to prove theorems~(\ref{T:bound1}) and (\ref{T:bound2}).

\section{Euclidean minima -- proof of the main results}
In this section we prove the main results of the paper, namely the upper bounds
for Euclidean minima stated in \S 3. Recall that for any number field $K$ of degree $n$,
we have $$ M(K) \leq \left(\frac{\tau_{\min}(\fok}{n}\right)^{\frac{n}{2}} \sqrt{D_K},$$
where $D_K$ is the absolute value of the discriminant of $K$. For $K \in \mathcal{A}$, we
now have an upper bound (see Corollary (\ref{c:tau-min})) and this will be used in the proofs.

\medskip
\noindent
{\it Proof of th. 3.1}
Set
\[
 f = \frac{p^{r+1}+p^r+1-e^2}{12 p^{r+1}}, \quad C = \sqrt{f}, \quad \varepsilon = \frac{rn}{\upsilon}.
\]
Then, by  the inequality~(\ref{E:gen-est}) and Corollary (\ref{c:tau-min}), we get
\[
 M(K) \leq \left(\frac{\tau(\fok)}{n}\right)^{\frac{n}{2}} \sqrt{D_K} \leq C^n \cdot (\sqrt{D_K})^{\varepsilon}.	
\]

First we shall prove that $\varepsilon(K)$ has the stated properties. If $r=1$, then $\upsilon = n-1$ and $\varepsilon = \frac{n}{n-1}$,
which implies that $\varepsilon \leq 2$ with the equality only for $n=2$. It is also clear that $\varepsilon \to 1$ if $n \to \infty$.
Assume that $r \geq 2$. A simple calculation shows that
\[
 2 \upsilon - rn = r d p^{r-1}- 2 d\left( \frac{p^{r-1}-1}{p-1}  \right) - 2.
\]
Clearly,
\[
 \frac{p^{r-1}-1}{p-1} \leq (r-1) p^{r-2},
\]
which gives
\[
 2 \upsilon - rn \geq  r d p^{r-1}- 2d (r-1) p^{r-2} - 2 = rdp^{r-2}(p-1) + p^{r-2}d-2 > 0,
\]
which implies that $\varepsilon < 2$. Another simple calculation shows that
\begin{equation}\label{E:rna}
 \varepsilon = \frac{r\left(1+\frac{1}{\upsilon}\right)}{\left( r - \delta \right)},
\end{equation}
where
\[
 \delta = \frac{p^{r-1}-1}{p^{r-1}(p-1)}.
\]
Observe that
$$
\ln n < r \ln p < r (p-1) \leq rp -1.
$$
Since $0 \leq \delta < \frac{1}{p}$, it follows that for $n \geq 3$ we have
$$
0 \leq \frac{r}{r-\delta} - 1< \frac{1}{rp -1} < \frac{1}{\ln n}.
$$
Thus we get
\begin{equation}\label{E:gt1}
\lim_{n \to \infty} \frac{r}{r - \delta} = 1.
\end{equation}
Finally, it follows from~(\ref{E:rna}) that $\upsilon + 1 = (r - \delta) n$.
Since $r - \delta > r - 1 \geq 1$, we get $\upsilon > n-1$. Consequently
\begin{equation}\label{E:gt2}
 \lim_{n \to \infty} \left( 1 + \frac{1}{\upsilon} \right) = 1.
\end{equation}
Combining the equalities~(\ref{E:gt1}), (\ref{E:gt2}), we obtain
\[
  \lim_{n_K \to \infty} \varepsilon(K) = 1.
\]

Now we prove the properties of $C(K)$. We have
 \begin{equation*}
  f =  \frac{p^{r+1}+p^r+1-e^2}{12 p^{r+1}} \leq \frac{p^{r+1}+p^r}{12 p^{r+1}} = \frac{p+1}{12p} \leq  \frac{1}{9}
 \end{equation*}
and hence $C \leq \frac{1}{3}$. We can also write
\[
   f =  \frac{p^{r+1}+p^r+1-e^2}{12 p^{r+1}} = \frac{1}{12} \left( 1+ \frac{1}{p} + \frac{1}{p^{r+1}} - \frac{e^2}{p^{r+1}} \right).
\]
If $r=1$ and $e$ is constant, then $f$ clearly approaches $\frac{1}{12}$ as $p \to \infty$ and hence $C(K) \to \frac{1}{2\sqrt{3}}$.  Assume now that $r \geq 2$. Since $1 \leq e \leq p-1$,
it follows that
\[
 0 < f -\frac{1}{12} \leq \frac{1}{12p}
\]
and thus
\[
 0 < C(K) - \frac{1}{2\sqrt{3}} \leq \frac{\sqrt{3}}{12p}.
\]
Consequently,
\[
 \lim_{p_K \to \infty} C(K) = \frac{1}{2\sqrt{3}}
\]
and this concludes the proof of the theorem.

\bigskip
\noindent
{\it Proof of th. 3.2.}
We shall use the same notation as in the proof of Theorem~(\ref{T:bound1}). In addition, we set
\[
 \omega(K)= C(K) \cdot (\sqrt{p})^{\delta + \frac{1}{n}}.
\]
A simple calculation using Corollary (\ref{c:tau-min}) and formulas derived in the proof of Theorem~(\ref{T:bound1}) gives
\[
 \left( \frac{\tau(\fok)}{n} \right)^{\frac{1}{2}} \leq \omega(K)
\]
Then, by the inequality~(\ref{E:gen-est}), we have
\[
  M(K) \leq \omega^n \sqrt{D_K}.
\]
If $r \geq 2$, then $\delta + \frac{1}{n} \leq \frac{2}{p}$ and hence
\[
1 < (\sqrt{p})^{\delta + \frac{1}{n}} \leq \sqrt[p]{p}.
\]
Consequently, using Theorem~(\ref{T:bound1}), we obtain
\begin{equation}\label{E:omega-lim}
\lim_{p \to \infty} \omega(K) = \lim_{p \to \infty} C(K) \cdot \lim_{p \to \infty} (\sqrt{p})^{\delta + \frac{1}{n}} = \frac{1}{2\sqrt{3}}.
\end{equation}
Moreover, using Theorem~(\ref{T:bound1}) and the fact that the sequence $\{\sqrt[p]{p}\}_{p \geq 3}$ is decreasing, we also get
\[
\omega(K) \leq C(K) \cdot \sqrt[p]{p} \leq \frac{1}{3} \sqrt[3]{3} = 3^{-2/3}.
\]
If $r=1$, then $\delta = 0$ and
\[
 (\sqrt{p})^{\delta + \frac{1}{n}}  = (\sqrt{p^e})^{\frac{1}{p-1}}.
\]
Thus assuming that $e$ is constant, we see that~(\ref{E:omega-lim}) holds as well.
This concludes the proof of the theorem.
\bigskip

\section{Abelian fields of prime conductor}

If $K$ is an abelian field of conductor $p^r$ with $r \ge 2$, then we have seen that
 $M(K) \leq 2^{-n} \sqrt{D_K}$, cf. (\ref{T:Minkowski}).
In particular, if $K$ is totally real, then Minkowski's conjecture holds for $K$. If $r = 1$, that
is if the conductor of $K$ is prime, then our results are less complete. The aim of this
section is to have a closer look at this case. As we will see, one can prove Minkowski's
conjecture in a number of special cases when $K$ is totally real.
\subsection{Totally real fields}
Let us consider  the set $\mathcal{S}_e$ of all totally real abelian fields of prime conductor such that $[\cyclo : K]=e$, where $e$ is an even positive integer.
The Dirichlet prime number theorem implies that the set $\mathcal{S}_e$ is infinite. By Theorem~(\ref{T:bound2}), we have
\[
 \lim_{p_K \to \infty} \omega(K) = \frac{1}{2\sqrt{3}}.
\]
In particular, for each $e$ there exists $N=N(e)$ such that for every field $K \in \mathcal{S}_e$ with $p_K > N$ we have
\[
 \omega(K) \leq \frac{1}{2}.
\]
and hence Minkowski's conjecture holds for these fields. The next result shows that we can take $N(e)=2e^2$.
\begin{proposition}\label{P:fixed-codim}
 Let $e$ be an even positive integer and $K \in \mathcal{S}_e$. If $p_K > 2e^2$, then
\[
 M(K) \leq 2^{-n} \sqrt{D_K}
\]
\end{proposition}

\noindent
{\it Proof.}
We shall use the same notation as in the proofs of theorems~(\ref{T:bound1}) and (\ref{T:bound2}). Additionally, let
\[
T = \set{(x,y) \in \br^2 \mid x \geq 2 \textnormal{ and } y \geq 2x+1}
\]
and  $h: T \to \br$ be a function given by
\[
 h(x,y) = \frac{y^2+y+1-x^2}{3 y^2} \cdot y^{\frac{x}{y-1}}.
\]
Then we have
\begin{equation*}
 \omega(K) = C(K) \cdot p^{1/2n} = \frac{1}{2} \cdot \sqrt{h(e,p)}.
\end{equation*}
Thus it is enough to show that $h(e,p) \leq 1$ for $p > 2e^2$. We set
\begin{align*}
 h_1(x,y)&=\frac{y^2+y+1-x^2}{3 y^2}\\
 h_2(x,y) &= y^{\frac{x}{y-1}}.
\end{align*}
For every $x \geq 2$ and every $y > 2(x^2-1)$ we have
\[
 \frac{\partial h_1}{\partial y}(x,y)= \frac{1}{3} \left( \frac{2(x^2-1)}{y^3} - \frac{1}{y^2} \right) < 0.
\]
Furthermore, for every $x \geq 2$ and every $y \geq 2 x + 1$ we have
\[
 \frac{\partial h_2}{\partial y}(x,y) = \frac{x h_2(x,y)}{(y-1)^2} \left( \frac{y-1}{y} - \ln y \right) < 0.
\]
Consequently, for every $x \geq 2$ and every $i=1,2$ the function $y \mapsto h_i(x,y)$ is positive and decreasing on the interval $[2x^2, \infty)$. Hence the function
$y \mapsto h(x,y)$ is decreasing on the interval $[2x^2, \infty)$. Moreover,
\[
 h(x,2x^2)  = (2x^2)^{\frac{x}{2x^2-1}} \cdot \frac{4x^4+x^2+1}{12x^4} < (2x^2)^{\frac{1}{2x-1}} \frac{1}{2} \leq 1.
\]
Consequently, $h(x,y) < 1$ for all $x \geq 2$ and $y \geq 2 x^2$. The result follows.
\bigskip

We can have $h(e,p) \leq 1$ even if $p < 2e^2$, which in many cases allows us to show that Minkowski's conjecture holds for every $K \in \mathcal{S}_e$.

\begin{example}
 If $e \leq 1202$ is an even integer, then Minkowski's conjecture holds for every $K \in \mathcal{S}_e$.
\end{example}

\noindent
{\it Proof.}
 If $p > 2 e^2$, then the result follows from Proposition~(\ref{P:fixed-codim}). For $p < 2 e^2$, it has been verified using Magma Computational Algebra System~\cite{magma} that either $h(e,p) \leq 1$ or $n_K \leq 8$.
In the first case, the result follows from the proof of Proposition~(\ref{P:fixed-codim}). In the second case, it follows from the fact that Minkowski's conjecture is known to hold for fields of degree not exceeding $8$.

\bigskip The previous results are based on upper bounds of  $\tau_{\rm min}(\fok)$ obtained
through the lattice $(\fok,q)$.
Another approach is to estimate $\tau_{\rm min}(\fok)$ using a scaling factor $\alpha$,
giving rise to the lattice $(\fok,q_{\alpha})$, see
\S4. A computation shows that if $p < 400$, then for an appropriate $\alpha \in \mathcal P$
the lattice $(\fok,q_{\alpha})$ is isomorphic to the unit lattice. Then \cite{bayer}, Corollary (5.5)
implies that Minkowski's conjecture holds.
\bigskip
\subsection{Totally imaginary fields}
If $\przystaje{p}{3}{4}$, then $K=\bq(\sqrt{-p}) \in \mathcal{A}$. Using formulas derived in the proof of Theorem~(\ref{T:bound1}), we get
\[
 M(K) \leq \frac{(p+1)^2}{16 p}.
\]
Note that this  bound is known to be the exact
value of the Euclidean minimum of $K$ (see for instance Proposition (4.2) in~\cite{lemmermeyer}). In particular, the inequality
$$
 M(K) \leq 2^{-n} \sqrt{D_K}
$$
does not hold in general for number fields that are not totally real. Just as in the totally real case, we have
\begin{equation}\label{E:tot-im-1}
 M(K) \leq  3^{-2n/3} \sqrt{D_K} < 2^{-n} \sqrt{D_K}
\end{equation}
for all totally imaginary fields $K \in \mathcal{A}$ with composite conductors. If the conductor of $K$ is prime and 
$n_K > 2$, then by Theorem~(\ref{T:bound1}) we have
\begin{equation}\label{E:tot-im-2}
 M(K) \leq 3^{-n} (\sqrt{D_K})^{\varepsilon}
\end{equation}
with $\varepsilon < 2$. Note that given the asymptotic behavior of the expressions $\varepsilon(K)$, $C(K)$, $\omega(K)$,
using the formulas derived in the proofs of Theorems~(\ref{T:bound1}) and (\ref{T:bound2}) directly will often lead to better bounds.


%
%

\end{document}